\documentclass{article}
\usepackage{amssymb, amsmath, amsfonts, amsthm, bm, amsbsy, graphicx}
\DeclareMathOperator{\Tr}{Tr} \DeclareMathOperator{\adj}{adj}

\newtheorem{thm}{Theorem}[section]
\newtheorem{cor}[thm]{Corollary}
\newtheorem{claim}[thm]{Claim}

\begin{document}
\date{}
\title{Aloff-Wallach Spaces: Volumes, Curvatures, Injectivity Radii}
\author{Mikhail Alyurov}
\maketitle

\section{Introduction} One of the central themes of
Riemannian geometry is the study of how local properties
(curvature) of a Riemannian manifold affect its global
(topological or metric) properties. The most famous example of
this is the classical Gauss-Bonnet Theorem. Many results relating
local and global properties are based on the injectivity radius
estimates. One example of this sort is Klingenberg's injectivity
radius estimate for the quarter-pinched \footnote{A manifold is
called $\delta$-pinched if its sectional curvatures lie between
two positive constants whose ratio is bounded by $\delta.$}
(compact simply connected Riemannian) manifolds as the main part
of the proof of the Sphere Theorem (\cite{K2}, also chapter 13 of
\cite{dC}).
\par Looking at the Klingenberg's injectivity radius estimate for the
quarter pinched manifolds, one would like to get an injectivity
radius estimate for $\delta$-pinched compact simply connected
Riemannian manifolds with any $\delta\in (0,1].$ Actually, the
problem only exists for odd dimensional manifolds since, two years
before proving the injectivity radius estimate for the quarter
pinched manifolds, Klingenberg showed, in \cite{K1}, that, for any
compact simply connected even dimensional manifold $M$ with
positive sectional curvature $K_M,$ the injectivity radius $i(M)$
satisfies
$$i(M)\geq\frac{\pi}{\sqrt{\max K_M}}.$$ The first instinct is to
try to get an estimate depending only on $\delta$ and the
dimension. This turns out to be impossible. In fact, Aloff-Wallach
spaces provide a counterexample to such an estimate.
\par In \cite{KS}, Klingenberg and Sakai conjectured that, if one
fixed a compact simply connected differential manifold $M$ and
then considered all possible $\delta$-pinched Riemannian
structures on $M,$ then one should be able to find a uniform lower
bound for the injectivity radii of the obtained Riemannian
manifolds. In the positively pinched case, finding a lower bound
on the injectivity radius is the same as finding a lower bound on
the volume. Therefore, the conjecture can be reformulated as: ``A
sequence of $\delta$-pinched Riemannian structures on a given
compact simply connected differential manifold can not collapse,''
where ``collapse'' means ``volume goes to zero.'' In this form,
the problem asks for application of methods of Gromov-Hausdorff
convergence. This approach (and in particular usage of the
N-structures introduced in \cite{CGF}) brought significant success
in proving the conjecture under different special assumptions. In
particular, in \cite{PRT}, it is proven that Klingenberg-Sakai
conjecture holds if, instead of considering all possible metrics,
one considers only metrics with bounded distance function, and
\cite{FR} contains a proof of the conjecture for the manifolds
satisfying special topological condition, namely that the second
Betti number is zero. As far as we know, the conjecture in its
general form is still open.
\par In this paper, we are going to focus on a particular example
of $\delta$-pinched manifolds, which may be interesting in its own
right. The topic of our study - Aloff-Wallach spaces, which were
first introduced in \cite{AW}, are the quotients of $SU(3)$ by
various images of $S^1.$ In \cite{H}, Huang showed that there is
an infinite family of uniformly pinched simply connected
topologically distinct Aloff-Wallach spaces and then used
Cheeger's Finiteness Theorem \cite{C} to conclude that this family
does not have a common lower injectivity bound.
\par The main results of this paper are two-sided volume estimates
for all Aloff-Wallach spaces [Theorem \ref{the_est}] and sharp
(sectional) curvature estimates for the Aloff-Wallach spaces from
the family mentioned in the last paragraph [Theorem
\ref{curv_thm}]. The estimation of the volumes uses
\emph{generalized Euler angles} on $SU(3),$ and the sectional
curvature bounds are obtained using \emph{modified curvature
operators} and the computational procedures given by P\"uttmann in
\cite{P}. As an application of these results, we obtain
injectivity radii estimates [Corollary \ref{ir}], which, in
particular, give a different proof of the Huang's result.
\par This paper is a part of the author's forthcoming Ph.D.
thesis in the Mathematics Department of Columbia University. The
author would like to thank his advisor, D.H. Phong, for numerous
helpful suggestions and constant encouragement.

\section{Description of the Spaces and\\ Statement of the Results}
\label{description} For each pair of integers $p$ and $q,$ we
define the subgroup $T(p,q)$ of $SU(3)$ by
\begin{equation}
T(p,q)= \left\{\left.
\begin{pmatrix}e^{2\pi ip\theta}&0&0\\
0&e^{2\pi iq\theta}&0\\
0&0&e^{-2\pi i(p+q)\theta}
\end{pmatrix}
\right|\theta\in\mathbb{R}\right\}.
\label{T_pq}
\end{equation}
If at least one of the numbers $p$ and $q$ is not zero, the
subgroup $T(p,q)$ is nontrivial, and the factor space
$$W(p,q)=SU(3)/T(p,q)$$ is called an Aloff-Wallach space.
\par It is shown in \cite{AW} that, if neither of $p,$ $q,$ and $p+q$
is zero, $W(p,q)$ can be equipped with a positively curved metric.
The positively curved metric on this space is obtained by
deforming the metric induced by the Killing form. The standard
Killing metric $k$ on $SU(3)$ is given by the formula
$$k(X,Y)=\frac{1}{2}\Tr(XY^*)$$
for $X,Y\in T_I(SU(3))$ and then extended by left invariance. This
metric induces a $SU(3)$-invariant metric on $W(p,q)$ in the
following way. We decompose $T_I(SU(3))=\mathfrak{su}(3)$ as
$$\mathfrak{su}(3)=\mathfrak{T}\oplus\mathfrak{T}^{\perp},$$ where
$\mathfrak{T}$ is the Lie algebra of $T(p,q)$ and
$\mathfrak{T}^\perp$ is its orthogonal complement in
$\mathfrak{su}(3)$ with respect to $k.$ Let
$$\pi:SU(3)\to W(p,q)(=SU(3)/T(p,q))$$ be the canonical
projection. The differential of the canonical projection at the
identity, $d\pi_I,$ gives an isomorphism of
$\mathfrak{T}^{\perp}=\mathfrak{su}(3)/\mathfrak{T}$ and
$T_{T(p,q)}(W(p,q)).$ Therefore, we shall have a scalar product on
$T_{T(p,q)}(W(p,q))$ once we have a scalar product on
$\mathfrak{T}^{\perp}.$ In order to be able to extend this product
by left invariance and obtain a $SU(3)$-invariant metric on
$W(p,q),$ the scalar product must be $Ad_{T(p,q)}$-invariant. One
obvious way to get such a scalar product on $\mathfrak{T}^{\perp}$
is to restrict the scalar product given by $k$ from
$\mathfrak{su}(3)$ to
$\mathfrak{T}^{\perp}\subset\mathfrak{su}(3).$ More generally,
supposing that there is an orthogonal, with respect to $k,$
$Ad_{T(p,q)}$-invariant decomposition
$\mathfrak{T}^{\perp}=V_1\oplus V_2,$ we can deform $k$ to obtain
a new $Ad_{T(p,q)}$-invariant scalar product $\tilde{k}$ on
$\mathfrak{T}^{\perp}$ as follows
$$\tilde{k}(X,Y)=a_1k(X_1,Y_1)+a_2k(X_2,Y_2),$$ where $a_1$ and
$a_2$ are positive constants, $X_1$ and $X_2$ are the projections
of $X$ on $V_1$ and $V_2$, and analogously for $Y.$
\par The construction of the aforementioned metric on $W(p,q)$ is based on
a particular choice of subspaces $V_1$ and $V_2$ and constants
$a_1$ and $a_2.$ First, we choose $V_1$ and $V_2.$ The choice is
made in the following way to ensure that the decomposition
satisfies certain conditions, called ``condition II'' in
\cite{AW}, which guarantee that $\tilde{k}$ (with appropriately
chosen constants $a_1$ and $a_2$) will induce an $SU(3)$-invariant
positively curved metric on $W(p,q).$ More precisely, it is shown
in \cite{AW} that, if the $Ad_{T(p,q)}$-invariant orthogonal
decomposition $\mathfrak{T}^{\perp}=V_1\oplus V_2$ satisfy:
\begin{enumerate}
\item $[V_1,V_2]\subset V_2,$ \item
$[V_1,V_1]\subset\mathfrak{T}\oplus V_1,$ \item $[V_2,V_2]\subset
\mathfrak{T}\oplus V_1,$ \item for any pair of linearly
independent vectors $x=x_1+x_2$ and $y=y_1+y_2,$ with $x_i,y_i\in
V_i,$ $[x,y]=0$ implies $[x_1,y_1]\neq 0,$
\end{enumerate}
then the metric $\tilde{k}(X,Y)=a_1k(X_1,Y_1)+a_2k(X_2,Y_2)$ as
above has positive curvature for $a_2=1$ and any $a_1\in(0,1).$ We
shall refer to the list above as \emph{condition II.}
\par In order to choose $V_1$ and $V_2,$ we start with the subgroup $U$ of $SU(3)$ given by
$$U=\left\{\left.\begin{pmatrix} g&0\\0&(\det g)^{-1}
\end{pmatrix}\right|g\in U(2)\right\}.$$
Note that this subgroup contains $T(p,q).$ The Lie algebra
$\mathfrak{u}$ of $U$ is given by
$$\mathfrak{u}=\left\{\left.\begin{pmatrix} u&0\\0&-\Tr
u\end{pmatrix}\right|u\in \mathfrak{u}(2)\right\}.$$ Let us point
out that $\mathfrak{T}\subset \mathfrak{u},$ which follows from
$T(p,q)\subset U,$ or could be seen directly from the fact that
$$\mathfrak{T}=\left\{\left.\begin{pmatrix}2\pi ip\theta&0&0\\0&2\pi
iq\theta&0\\0&0&-2\pi
i(p+q)\theta\end{pmatrix}\right|\theta\in\mathbb{R}\right\} .$$ We
form the decomposition of $\mathfrak{T}^{\perp}$ by taking
\begin{equation}
\begin{split}
V_1&=\mathfrak{T}^{\perp}\cap\mathfrak{u},\\
V_2&=\mathfrak{u}^{\perp},
\end{split}
\label{dec_def}
\end{equation}
where $\mathfrak{u}^{\perp}$ is the orthogonal complement of
$\mathfrak{u}$ with respect to the Killing form $k.$ The fact that
$\mathfrak{T}^{\perp}=V_1\oplus V_2$ follows from the fact that
$\mathfrak{T}\subset \mathfrak{u}.$ A series of matrix
computations shows that $V_1$ and $V_2$ given by \eqref{dec_def}
are $Ad_T(p,q)$-invariant, and, if $pq>0,$ the decomposition
$\mathfrak{T}^{\perp}=V_1\oplus V_2$ (with $V_1$ and $V_2$ given
by \eqref{dec_def}) satisfy condition II.
\par We complete the construction of our particular version of the
positively curved metric on $W(p,q)$ by picking $a_1=1/2$ and
$a_2=1,$ which makes
$$\tilde{k}(X,Y)=\frac{1}{2}k(X_1,Y_1)+k(X_2,Y_2).$$
\par Since, for $pq>0,$ the decomposition
$\mathfrak{T}^{\perp}=V_1\oplus V_2,$ with $V_1$ and $V_2$ given
by \eqref{dec_def}, satisfies condition II, Theorem 2.4 of
\cite{AW} says that $SU(3)$-invariant metric induced on $W(p,q)$
[with $pq>0$] by $\tilde{k}$ is positively curved; using this,
Theorem 3.2 of \cite{AW} shows that the result holds as long as
neither of $p,$ $q,$ or $p+q$ is zero.
\par In \cite{H}, Huang proved that the curvature of $W(p,q)$ (with
this metric) depends only on the ratio $p/q$ and established that
the curvature of $W(1,1)$ is pinched between $2/37$ and $29/8.$
Using this, he showed that the Aloff-Wallach spaces $W(i,i+1),$
with $i$ sufficiently big, are uniformly pinched, simply
connected, and topologically distinct, and, therefore, do not have
a common lower injectivity radius.
\par Let us now formulate our results precisely.
\begin{thm}[Volumes] \label{the_est} If $Vol(W(p,q))$ denotes the volume
of the Aloff-Wallach space $W(p,q)$ with respect to the metric
that we have just chosen, then
$$\frac{\sqrt{3}\pi^4\gcd(p,q)}{32\sqrt{p^2+q^2+pq}}\leq Vol(W(p,q))
\leq\frac{\sqrt{3}\pi^4\gcd(p,q)}{2\sqrt{p^2+q^2+pq}}.$$
\end{thm}
\begin{thm}[Curvatures]
\label{curv_thm} For any positive integer $n,$ the sectional
curvature of the Aloff-Wallach space $W(n,n+1),$ satisfies the
sharp inequality $$c(n)\leq K(W(n,n+1))\leq C(n),$$ where
\begin{align*}
c(n)&=\frac{17+63n+63n^2}{16+48n+48n^2}\\
&-\frac{1}{16}\left[\frac{(7+33n+33n^2)^2}{(1+3n+3n^2)^2}+4\left(\frac{9(1+2n)}{\sqrt{3+9n+9n^2}}\right.\right.\\
&+\left\{(32+552n+3132n^2+8037n^3+9648n^4+4401n^5)\sqrt{3+9n+9n^2}\right.\\
&-\sqrt{3n}(16+60n+57n^2)(-56-555n-1935n^2-1620n^3+7173n^4\\
&\hphantom{-\sqrt{3n}(16+60n+57n^2)(}\left.\left.+22788n^5+26649n^6+11907n^7)^{1/2}
\vphantom{\sqrt{9n^2}}\right\}\right/\\
&\hphantom{+}\left\{64+672n+2916n^2+6624n^3+8181n^4+4995n^5+999n^6\right\}
\left.\vphantom{\frac{9(1+2n)}{\sqrt{3+9n+9n^2}}}\right)^2
\left.\vphantom{\frac{(7+33n+33n^2)^2}{(1+3n+3n^2)^2}}\right]^{1/2}
\end{align*}
and $$C(n)=4-\frac{9n^2}{8(1+3n+3n^2)},$$ which implies non-sharp
inequality
$$\frac{1}{25}\leq K(W(n,n+1))\leq 4.$$
\end{thm}
For any Riemannian manifold $M,$ let $i(M)$ denote its injectivity
radius.
\begin{cor}[Injectivity Radii]
\label{ir} The injectivity radii of the various Aloff-Wallach
spaces satisfy the following inequalities:
\begin{enumerate}
\item $i(W(1,1))\geq 4.65\cdot10^{-5};$ \label{ir_1} \item
$i(W(n,n+1))\geq\dfrac{3\sqrt{3}\pi(c(n))^3}{32\sqrt{3n^2+3n+1}},$\\
where $c(n)$ is the functions from Theorem \ref{curv_thm};
\label{ir_n} \item $i(W(p,q))\leq
\pi\cdot\left[\dfrac{3\sqrt{3}\gcd(p,q)}{2\sqrt{p^2+q^2+pq}}\right]^{1/7}.$
\label{ir_ab}
\end{enumerate}
\end{cor}
Theorem \ref{the_est} is established in section \ref{vol_sec},
Theorem \ref{curv_thm} in section \ref{curv_sec}, and Corollary
\ref{ir} in section \ref{ir_sec}.

\section{Volume of $W(p,q)$}
\label{vol_sec}
\subsection{Preliminary considerations}
In order estimate the volume of $W(p,q)$ and prove Theorem
\ref{the_est}, we are going to use the following result: If
$\pi:(G,g)\to (M,f)$ is a Riemannian submersion, then
$$Vol(G,g)=\int_M{Vol\left(\pi^{-1}(x)\right)\sqrt{\det f(x)}\,dx},$$
which is given as Corollary II.5.7 in \cite{S}. If in addition G
is a Lie group and $M$ is its homogeneous space, say $M=G/H,$ then
points of $M$ are left cosets: for any $x\in M$ there exists $g\in
G$ such that $x=[gH],$ and $\pi^{-1}(x)=\pi^{-1}([gH])=gH.$ If
further the metric on $G$ is left-invariant, all $gH$ mentioned in
the previous sentence are isometric, and, in particular, their
volumes are equal. Thus, if $\pi:G\to G/H$ is the canonical
projection,
\begin{equation}
Vol(G,g)=\int_{G/H}{Vol(H)\sqrt{\det f(x)}\,dx}=Vol(H)\cdot
Vol(G/H,f). \label{subm_vol}
\end{equation}
\par In order to apply this formula to $(G,H)=(SU(3),T(p,q)),$
we need to pick a metric on $SU(3)$ such that the canonical
projection $$\pi:SU(3)\to W(p,q)$$ is a Riemannian submersion.
Such a metric is induced by the scalar product
$$w:\mathfrak{su}(3)\times\mathfrak{su}(3)\to\mathbb{R}$$ defined
by
$$w(X,Y)=k(X_{\mathfrak{T}},Y_{\mathfrak{T}})+\frac{1}{2}k(X_1,Y_1)+k(X_2,Y_2),$$
where $X_{\mathfrak{T}}\in\mathfrak{T},$ $X_1\in V_1,$ $X_2\in
V_2,$ and analogously for $Y.$ The actual metric, which we are
also going to call $w,$ is given by extending this scalar product
by left invariance.
\par{\footnotesize The fact that $\pi:(SU(3),w)\to (W(p,q),\tilde{k})$ is a Riemannian
submersion at $I\in SU(3)$ follows from the definitions of $w$ and
$\tilde{k}:$ With the identification of $T_{T(p,q)}(W(p,q))$ and
$\mathfrak{T}^{\perp}$ that was made in order to construct the
metric on $W(p,q),$ $d\pi_I:(\mathfrak{su}(3),w)\to
(T_{T(p,q)}(W(p,q)),\tilde{k})$ is an orthogonal projection. Since
the metrics on $SU(3)$ and $W(p,q)$ are left-invariant, this
implies that $\pi$ is a Riemannian submersion everywhere.}
\par The application of \eqref{subm_vol} to $\pi:(SU(3),w)\to
(W(p,q),\tilde{k})$ yields
\begin{equation}
Vol(W(p,q),\tilde{k})=Vol(SU(3),w)/Vol(T(p,q)). \label{wal vol
frac}
\end{equation}
\par Now, our goal is to estimate the volume of
$SU(3)$ (in metric $w$) and compute the length of $T(p,q).$ In
view of \eqref{T_pq}, the tangent vector
$\mathfrak{t}=\mathfrak{t}(\theta)$ of $T(p,q)$ is given by
$$\mathfrak{t}=
\begin{pmatrix}
2\pi ip e^{2\pi ip\theta}&0&0\\
0&2\pi iq e^{2\pi iq\theta}&0\\
0&0&-2\pi i(p+q) e^{-2\pi i(p+q)\theta}
\end{pmatrix}.$$
Therefore,
\begin{equation}
Vol(T(p,q))=\int_0^{1/\gcd(p,q)}\sqrt{k(\mathfrak{t},\mathfrak{t})}\,d\theta
=\frac{2\pi}{\gcd(p,q)}\sqrt{p^2+q^2+pq}.
\label{orbit length}
\end{equation}
\subsection{Volume of $SU(3)$}
\subsubsection{Euler angle parametrization}
In order to compute the volume of $SU(3),$ we are going to
introduce the generalized Euler angles. Before we start describing
the parametrization of $SU(3),$ let us recall the Euler angles on
$SU(2).$ In the case of $SU(2),$ one uses the Pauli matrices
$\sigma_i$ given by
\begin{align*}
&\sigma_1=\begin{pmatrix}0&1\\1&0\end{pmatrix},&
&\sigma_2=\begin{pmatrix}0&-i\\i&0\end{pmatrix},&
&\sigma_3=\begin{pmatrix}1&0\\0&-1\end{pmatrix}
\end{align*}
to write a generic $s\in SU(2)$ as
$$s=s(\phi,\theta,\psi)=e^{\frac{i}{2}\phi\sigma_3}e^{\frac{i}{2}\theta\sigma_2}e^{\frac{i}{2}\psi\sigma_3},$$
thus parameterizing $SU(2)$ (outside a set of measure zero) by
$\phi,$ $\theta,$ and $\psi.$ Let us recall how the coordinate
ranges of this parametrization are found.
\par It follows directly from the definition of $SU(2)$ that
$$SU(2)=\left\{\begin{pmatrix}a&b\\-\bar b&\bar a\end{pmatrix}:\left\{a,b\in\mathbb{C}\ \&\
|a|^2+|b|^2=1\right\}\right\}.$$ On the other hand, since
$$
e^{\frac{i}{2}\theta\sigma_2}=\begin{pmatrix}\cos\frac{\theta}{2}&\sin\frac{\theta}{2}\\
-\sin\frac{\theta}{2}&\cos\frac{\theta}{2}\end{pmatrix},
$$
\begin{align*}
s(\phi,\theta,\psi)&=e^{\frac{i}{2}\phi\sigma_3}e^{\frac{i}{2}\theta\sigma_2}e^{\frac{i}{2}\psi\sigma_3}\\
&=\begin{pmatrix}e^{\frac{i}{2}(\phi+\psi)}\cos\frac{\theta}{2}&e^{\frac{i}{2}(\phi-\psi)}\sin\frac{\theta}{2}\\
-e^{-\frac{i}{2}(\phi-\psi)}\sin\frac{\theta}{2}&e^{-\frac{i}{2}(\phi+\psi)}\cos\frac{\theta}{2}\end{pmatrix}.
\end{align*}
Therefore, in order to find the ranges of the Euler angles, we
need to find three intervals $I_{\phi},$ $I_{\theta},$ and
$I_{\psi}$ so that, outside a set of measure 0, there is a
diffeomorphism
$$f:\{(a,b)\in\mathbb{C}^2:|a|^2+|b|^2=1\}\to I_{\phi}\times
I_{\theta}\times I_{\psi},$$ with the property that, if
$f(a,b)=(\phi, \theta, \psi),$ then
\begin{equation*}
\begin{split}
a&=e^{\frac{i}{2}(\phi+\psi)}\cos\frac{\theta}{2},\\
b&=e^{\frac{i}{2}(\phi-\psi)}\sin\frac{\theta}{2}.
\end{split}
\label{Euler angles}
\end{equation*}
\par Writing $(a,b)$ in the form $(|a|e^{i\alpha},|b|e^{i\beta}),$
with $\alpha,\beta\in[0,2\pi),$ we define
$f:(a,b)\mapsto(f_\phi(a,b), f_\theta(a,b), f_\psi(a,b))$ by
$f_\theta(a,b)=2\arccos(|a|)$ and
$$
(f_\phi(a,b),f_\psi(a,b))=
\begin{cases}
(\alpha+\beta, \alpha-\beta)&\text{if $\alpha\geq\beta$},\\
(\alpha+\beta-2\pi, \alpha-\beta+2\pi)&\text{if $\alpha<\beta$ and $\alpha+\beta\geq2\pi$},\\
(\alpha+\beta+2\pi, \alpha-\beta+2\pi)&\text{if $\alpha<\beta$ and $\alpha+\beta<2\pi$},\\
\end{cases}
$$
which gives
$$I_{\phi}\times
I_{\theta}\times I_{\psi}=[0,4\pi)\times[0,\pi)\times[0,2\pi).$$
\par In the case of $SU(3),$ Gell-Mann matrices $\lambda_i$ are used in
place of Pauli matrices. The Gell-Mann matrices that are used in
the parametrization are
\begin{align*}
&\lambda_2=\begin{pmatrix}0&-i&0\\i&0&0\\0&0&0\end{pmatrix},&
&\lambda_3=\begin{pmatrix}1&0&0\\0&-1&0\\0&0&0\end{pmatrix},\\
&\lambda_5=\begin{pmatrix}0&0&-i\\0&0&0\\i&0&0\end{pmatrix},&
&\lambda_8=\frac{1}{\sqrt{3}}\begin{pmatrix}1&0&0\\0&1&0\\0&0&-2\end{pmatrix}.
\end{align*}
We claim that any $g\in SU(3)$ can be written as
$$g=g(\phi,\theta,\psi,\xi,\alpha,\beta,\gamma,\tau)=
s(\phi,\theta,\psi)e^{i\lambda_5\xi}
s(\alpha,\beta,\gamma)e^{\frac{i\sqrt{3}}{2}\lambda_8\tau},$$
where
$$s(x,y,z)=e^{\frac{i}{2}x\lambda_3}
e^{\frac{i}{2}y\lambda_2}e^{\frac{i}{2}z\lambda_3}$$ is the
Euler-angle parametrization of $SU(2)\subset SU(3).$ The
coordinates $\phi,$ $\theta,$ $\psi,$ $\xi,$ $\alpha,$ $\beta,$
$\gamma$ and, $\tau$ as above are called generalized Euler angles.
\par
Direct computation, which can be done painlessly with the help of
\emph{Mathematica,} shows that, for any choice of parameters,
$$g(\phi,\theta,\psi,\xi,\alpha,\beta,\gamma,\tau)\in SU(3).$$ To
find the ranges of the coordinates we look at the matrix elements
of $$g=s(\phi,\theta,\psi)e^{i\lambda_5\xi}
s(\alpha,\beta,\gamma)e^{\frac{i\sqrt{3}}{2}\lambda_8\tau},$$ and,
by considerations similar to the ones used in the $SU(2)$ case, we
establish that
\begin{align*}
\beta,\theta,\xi&\in[0,\pi);&\alpha,\phi&\in[0,4\pi);&\gamma,\psi&\in[0,2\pi).
\end{align*}
\subsubsection{Estimation}
\par In computation of the volume, we are going to use the coordinate
vector fields, which are obtained by differentiating
$g=g(\phi,\theta,\psi,\xi,\alpha,\beta,\gamma,\tau),$ for example
$\partial_{\alpha}=\frac{\partial g}{\partial\alpha}.$ Using these
as a basis, we, \emph{theoretically}, could compute the
determinant of the metric $w$ and then integrate the square root
of this determinant to get the volume. In reality, however, this
computation is too complicated even for \emph{Mathematica} to
handle. Therefore, we shall settle for the estimate of the volume
near the volume in the Killing metric, whose volume element turns
out to be given by a nice formula. Recall that we decomposed the
Lie algebra of $SU(3)$ as
$$\mathfrak{su}(3)=\mathfrak{T}\oplus V_1\oplus V_2,$$ and defined the
Wallach scalar product $w$ on $\mathfrak{su}(3)$ by
$$w(X,Y)=k(X_{\mathfrak{T}},Y_{\mathfrak{T}})+\frac{1}{2}k(X_1,Y_1)+k(X_2,Y_2),$$
where $k$ is the killing form on $SU(3)$ and the subscripts denote
the projections on the corresponding subspaces. Note that $w$ can
be easily estimated in terms of $k:$
$$\frac{1}{2}k(X,X)\leq w(X,X)\leq k(X,X),$$
and, hence,
\begin{equation}
\frac{1}{16}Vol(SU(3),k)\leq Vol(SU(3),w)\leq Vol(SU(3),k).
\label{vol est}
\end{equation}
Therefore, once we know the volume of $SU(3)$ in the Killing
metric, we'll have a two-sided estimate on the volume in the
Wallach metric, which is our main goal.
\par Using the strategy described in the beginning of the last paragraph,
we compute, with the help of \emph{Mathematica,} that the volume
element of the Killing metric, $k,$ at a generic point is given by
the formula
$$dV=\frac{\sqrt{3}}{512}\sin\beta\sin\theta\sin\xi\sin^2\frac{\xi}{2}.$$
Integrating this formula over the ranges of the generalized Euler
angles, we get
$$Vol(SU(3),k)=\sqrt{3}\pi^5.$$
\par Inserting this into \eqref{vol est}, we get the following
two-sided estimate for the volume of $SU(3)$ with the
Aloff-Wallach metric:
$$\frac{\sqrt{3}}{16}\pi^5\leq Vol(SU(3),w)\leq\sqrt{3}\pi^5.$$
Combining this with \eqref{wal vol frac} and \eqref{orbit length},
we get the promised two-sided  estimate for the volume of the
Aloff-Wallach spaces:
\begin{equation}
\frac{\sqrt{3}\pi^4\gcd(p,q)}{32\sqrt{p^2+q^2+pq}}\leq Vol(W(p,q))
\leq\frac{\sqrt{3}\pi^4\gcd(p,q)}{2\sqrt{p^2+q^2+pq}},
\label{main_est}
\end{equation}
completing the proof of Theorem \ref{the_est}.

\section{The Pinching of W(n,n+1)}
\label{curv_sec}
\subsection{General Remarks about Curvatures}
Our estimation of the curvature (and computation of the
pinching)is based on the procedure given in \cite{P}. The central
tool of the procedure is \emph{modified curvature operators.} Let
us describe the relevance of the modified curvature operators to
the estimation of the sectional curvature.
\par For any Riemannian manifold $(M,g)$ and the corresponding
Levi-Civita connection $\nabla,$ one defines the Riemann curvature
tensor
$Rm:\Gamma(TM)\times\Gamma(TM)\times\Gamma(TM)\times\Gamma(TM)\to\mathbb(R)$
by the formula
$$Rm(X,Y,Z,W)=g(\nabla_X\nabla_{Y}Z-\nabla_Y\nabla_{X}Z-\nabla_{[X,Y]}Z,W).$$
It can be shown that $Rm$ can be used to define a symmetric
bilinear form on the bundle of bivectors,
$\widehat{Rm}:\Lambda_{2}TM\otimes\Lambda_{2}TM\to\mathbb{R},$ by
$$\widehat{Rm}(X\wedge Y,W\wedge Z)=Rm(X,Y,Z,W).$$ The
self-adjoint linear operator associated to this symmetric bilinear
form is called the \emph{curvature operator} and will be denoted
by $\mathfrak{R}.$ In other words, $\mathfrak{R}$ is defined by
the equality
$$\hat{g}(\mathfrak{R}(X\wedge Y),W\wedge Z)=\widehat{Rm}(X\wedge Y,W\wedge
Z),$$ where $\hat{g}$ is the metric induced by $g$ on
$\Lambda_{2}TM.$
\par The \emph{sectional curvature} $K$ of a given Riemannian manifold
$(M,g)$ is a function that associates to any pair of linearly
independent vectors $\{X,Y\}\subset T_pM$ (for some $p\in M$) a
number
$$K(X,Y)=\frac{\widehat{Rm}(X\wedge Y,X\wedge
Y)}{\hat{g}(X\wedge Y,X\wedge Y)}.$$ The value of the sectional
curvature depends only on the 2-plane spanned by $X$ and $Y$
(which makes it possible to write $K(X\wedge Y)$ in place of
$K(X,Y)$). Therefore, in order to estimate the sectional
curvature, it is enough to look at its values on the orthonormal
pairs of vectors. For an orthonormal pair of vectors $\{X,Y\}$,
$$K(X,Y)=\widehat{Rm}(X\wedge Y,X\wedge Y)=\hat{g}(\mathfrak{R}(X\wedge Y),X\wedge Y),$$
which lies between the smallest and the biggest eigenvalues of
$\mathfrak{R}.$ Thus, one way to estimate the sectional curvature
is to compute the eigenvalues of the curvature operator.
Unfortunately, this estimate is not optimal because an eigenvector
of $\mathfrak{R}$ might happen to be a bivector that can not be
written as a wedge of two tangent vectors. In particular, the
smallest eigenvalue of the curvature operator on $W(1,1)$ is
negative.
\subsection{Modified Curvature Operators}
\label{mod_curv_op_ss} The shortcomings of the curvature operator
method of estimating the sectional curvature described at the end
of the last subsection can be overcome if one considers
\emph{modified curvature operators} in place of the curvature
operator. The construction of the modified curvature operators, is
based on the function $i:\Lambda^{4}TM\to S^{2}(\Lambda_{2}TM)$
that assigns to each 4-form $\Omega$ a symmetric bilinear form (on
the space of bivectors) $i(\Omega)$ defined by
$[i(\Omega)](\alpha_1,\alpha_2)=\Omega(\alpha_1\wedge\alpha_2).$
Now, for each $\Omega\in\Lambda^{4}TM,$ we define modified Riemann
curvature tensor by $Rm_\Omega=Rm+i(\Omega).$ The self-adjoint
linear operator associated to the symmetric bilinear form
$Rm_\Omega$ is called a modified curvature operator and is denoted
by $\mathfrak{R}_\Omega.$ Since $Rm_\Omega(X\wedge Y,X\wedge
Y)=Rm(X\wedge Y,X\wedge Y)$ for all $\Omega\in\Lambda^{4}TM,$ the
sectional curvature is controlled by the eigenvalues of the
modified curvature operators:
\begin{equation}
\lambda_{min}(\mathfrak{R}_{\Omega_1})\leq
K\leq\lambda_{max}(\mathfrak{R}_{\Omega_2}),
\label{cop_est_sec}
\end{equation}
where $\Omega_1$ and $\Omega_2$ are any two 4-forms.
\par Let us describe how this inequality can be used to estimate
the sectional curvature of the Aloff-Wallach spaces. The strategy
is to estimate the curvature on certain subspaces of
$\Lambda_{2}(T_{T(p,q)}W(p,q))$ and then to show that the bounds
are stricter then the bounds given by the corresponding
eigenvalues of some modified curvature operators. To formulate
this more precisely, we introduce
\begin{align*}
\hat{\lambda}&:=\inf\{K(\omega)|\omega\in G\cap E_1\},\\
\bar{\lambda}&:=\inf\{K(\omega)|\omega\in G\cap E_2\},\\
\Lambda_j&:=\inf\{K(\omega)|\omega\in G\cap F_i\},
\end{align*}
where $j\in\{0,1,2\},$ $G$ is the Grassmannian of oriented
2-planes, and $E_i$ and $F_i$ are subspaces of
$\Lambda_{2}(T_{T(p,q)}W(p,q)).$ It is established in \cite{P}
that $\min\{\hat{\lambda},\bar{\lambda}\}$ is weakly smaller then
the minimal eigenvalue of certain modified curvature operator.
Using inequality \eqref{cop_est_sec}, this implies that
$\min\{\hat{\lambda},\bar{\lambda}\}$ bounds the sectional
curvature from below, and, since it is clear from the definition
that $\hat{\lambda}$ and $\bar{\lambda}$ are weakly larger than
the minimum of the (unrestricted) sectional curvature, $K_{min},$
we get $\min\{\hat{\lambda},\bar{\lambda}\}=K_{min}.$ Similar
reasoning works for $\Lambda_j$ and $K_{max}.$
\par Sections 5 of \cite{P} gives concrete recipes for computing
$\hat{\lambda},$ $\bar{\lambda},$ and $\Lambda_j.$ In his paper,
P\"{u}tmann uses these recipes to determine the optimal pinching
among certain class of metrics on the Aloff-Wallach spaces. We
shall employ these procedures to determine the minimum and the
maximum of the sectional curvatures on the Aloff-Wallach spaces
with the metric $\tilde{k}$, which we defined in section
\ref{description}.
\par In order to proceed to the computations,
we need to introduce some notation, which we are going to take
from \cite{P}, but adapt to our case. In proposition 4.10 (of
\cite{P}), P\"{u}ttmann introduces quantities $a_j,$ $b_j,$ $c_j,$
$d_j,$ and $\xi_j,$ with $j\in\{0,1,2\}$ (and shows that they are
the matrix elements of the curvature operators restricted to
various subspaces of $\Lambda_{2}TM$). Here are the definitions of
these quantities for an Aloff-Wallach space $W(p,q)$ adapted to
our choice of metric:
\begin{align*}
a_0&=8-\frac{9(p+q)^2}{2(p^2+pq+q^2)},&
a_1&=4-\frac{9p^2}{8(p^2+pq+q^2)},&
a_2&=4-\frac{9q^2}{8(p^2+pq+q^2)};\\
b_0&=-2-\frac{9pq}{8(p^2+pq+q^2)},&
b_1&=-\frac{10p^2+pq+q^2}{4(p^2+pq+q^2)},&
b_2&=-\frac{p^2+pq+10q^2}{4(p^2+pq+q^2)};\\
c_0&=\frac{3(p+q)^2}{2(p^2+pq+q^2)},&
c_1&=\frac{3p^2}{8(p^2+pq+q^2)},&
c_2&=\frac{3q^2}{8(p^2+pq+q^2)};\\
d_0&=\frac{5}{8},& d_1&=\frac{1}{8},&
d_2&=\frac{1}{8};\\
\xi_0&=-\frac{3\sqrt{3}(p+q)}{8\sqrt{p^2+pq+q^2}},&
\xi_1&=\frac{\sqrt{3}(2p+q)}{8\sqrt{p^2+pq+q^2}},&
\xi_2&=\frac{\sqrt{3}(p+2q)}{8\sqrt{p^2+pq+q^2}}.\\
\end{align*}
\par Here are more details on these quantities: It turns out that it
is enough to consider \emph{invariantly} modified curvature
operators, which are the operators $\mathfrak{R}_\Omega$ with
$\hat{T}^2$-invariant $\Omega,$ where $\hat{T}^2$ is the extension
of $T^2\subset SU(3)$ by the complex conjugation on $SU(3).$ To
decompose these operators, their domain
$\Lambda_{2}(T_{T(p,q)}W(p,q))=\Lambda_{2}\mathfrak{T}^{\perp}$,
is decomposed into the sum of $\hat{T}^2$-invariant subspaces by
first identifying $\mathfrak{T}^{\perp}$ with
$\mathbb{R}\oplus\mathbb{C}^3,$ and then decomposing the
corresponding space of bivectors as
$$\Lambda_2(\mathfrak{T}^{\perp})=
\mathbb{R}^3\oplus\left(\oplus_{j=0}^{2}V_j\right)
\oplus\left(\oplus_{j=0}^{2}\left(\mathbb{C}_j^a\oplus\mathbb{C}_j^b\right)\right),$$
where the only non-obvious terms $V_j$ are copies of $\mathbb{C}.$
Then it is shown that $\hat{T}^2$-invariant 4-forms on
$\mathfrak{T}^{\perp}$ are parameterized by four real numbers.
Calling these numbers $\eta_0,$ $\eta_1,$ $\eta_2,$ and $\xi$ and
writing $\mathfrak{R}(\bm{\eta},\xi)$ for
$\mathfrak{R}_{\Omega(\eta_0,\eta_1,\eta_2,\xi)},$ one gets the
following decomposition for the invariantly modified curvature
operator:
\begin{align*}
\mathfrak{R}(\bm{\eta},\xi)|_{V_j}&=\eta_j,\\
\mathfrak{R}(\bm{\eta},\xi)|_{\mathbb{R}^3}&=\begin{pmatrix}a_0&b_2-\eta_2&b_1-\eta_1\\
b_2-\eta_2&a_1&b_0-\eta_0\\b_1-\eta_1&b_0-\eta_0&a_2\end{pmatrix},\\
\mathfrak{R}(\bm{\eta},\xi)|_{\mathbb{C}_j^a\oplus\mathbb{C}_j^b}&=\begin{pmatrix}c_j&\sqrt{2}(\xi_j-\xi)\\
\sqrt{2}(\xi_j-\xi)&2d_j-\eta_j\end{pmatrix}.\\
\end{align*}
We can also mention that
\begin{align*}
E_1&=\left(\oplus_{j=0}^{2}V_j\right)
\oplus\left(\oplus_{j=0}^{2}\left(\mathbb{C}_j^a\oplus\mathbb{C}_j^b\right)\right),\\
E_2&=\mathbb{R}^3\oplus\left(\oplus_{j=0}^{2}V_j\right),\\
F_j&=\mathbb{R}^3\oplus\mathbb{C}_j^a\oplus\mathbb{C}_j^b.
\end{align*}
\subsection{Minimal Curvature}
\label{ss_min_curv} As we explained in the previous subsection, in
order to find the minimum of the sectional curvature, we need to
compute $\hat{\lambda}$ and $\bar{\lambda}.$ Let us describe how
these two numbers are computed in \cite{P}. In order to compute
$\hat{\lambda},$ one considers three functions $\lambda_j(x)$
(where $j\in\{0,1,2\}$), each of which is the smallest root of the
corresponding polynomial
$$P_{x}(\lambda)=\det(\mathfrak{R}(\lambda,\xi)|_{\mathbb{C}_j^a\oplus\mathbb{C}_j^b}-\lambda
I),$$ where $\mathfrak{R}(\lambda,\xi)$ is the modification with
$(\eta_0,\eta_1,\eta_2)=(\lambda,\lambda,\lambda).$ In terms of
these functions, $\hat{\lambda}$ can be computed as
$\hat{\lambda}=\max_{x}\min_{j}\{\lambda_j(x)\}.$ It follows from
the decomposition of the modified curvature operators shown in the
previous subsection that
$$\lambda_j(x)=\frac{c_j+d_j}{2}-\sqrt{\left(\frac{c_j-d_j}{2}\right)^2+(\xi_j-x)^2}.$$
Since $y=\lambda_j(x)$ are lower branches of hyperbolas with
maxima at $x=\xi_j,$ $\hat{\lambda}$ is achieved either at $\xi_j$
or at an intersection of two curves between their maxima.
\par We shall now compute $\hat{\lambda}$ for $W(n,n+1)$ (with n a positive integer).
First let us look at the maxima of $\lambda_j,$ that is
$\lambda_j(\xi_j).$ In the following table, we write $\xi_j(n)$
for $\xi_j(p,q)$ with $p=n$ and $q=n+1$ and $\lambda_j(\xi_j)$ for
$\lambda_j(\xi_j(n)).$
\begin{align*}
\xi_0(n)&=-\frac{9(1+2n)}{8\sqrt{3+9n+9n^2}},&
\xi_1(n)&=\frac{3(1+3n)}{8\sqrt{3+9n+9n^2}},&
\xi_2(n)&=\frac{3(2+3n)}{8\sqrt{3+9n+9n^2}};\\
\lambda_0(\xi_0)&=\frac{5}{8},&
\lambda_1(\xi_1)&=\frac{3n^2}{8+24n+24n^2},&
\lambda_2(\xi_2)&=\frac{1}{8}.
\end{align*}
The numbers from the second row can not be $\hat{\lambda}$ since,
for any $j\in\{0,1,2\},$ there exists $i\in\{0,1,2\},$ with $i\neq
j,$ such that $\lambda_j(\xi_j)>\lambda_i(\xi_j).$ Namely,
$$\lambda_1(\xi_0)=\frac{1+3n+6n^2-\sqrt{193+1446n+4149n^2+5508n^3+2916n^4}}{16+48n+48n^2}<0,$$
which implies $\lambda_1(\xi_0)<\lambda_0(\xi_0)$ and
$\lambda_0(\xi_0)\neq\hat{\lambda};$
\begin{align*}
\lambda_0(\xi_1)&=\frac{17+63n+63n^2-\sqrt{241+1902n+5691n^2+7686n^3+4005n^4}}{16+48n+48n^2}\\
&<\frac{3n^2}{8+24n+24n^2}=\lambda_1(\xi_1)
\end{align*}
for any $n\geq 1$ which implies
$\lambda_1(\xi_1)\neq\hat{\lambda}$ by the same logic as above;
and, finally,
$$\lambda_0(\xi_1)=\frac{17+63n+63n^2-\sqrt{349+2442n+6663n^2+8334n^3+4005n^4}}{16+48n+48n^2}<0$$
shows $\lambda_2(\xi_2)\neq\hat{\lambda}.$
\par The other candidates for $\hat{\lambda}$ are the intersections, so let us
consider these. Setting $\lambda_0(x)=\lambda_1(x),$ we get two
roots:
\begin{align*}
x^{01}_1&=[(32+552n+3132n^2+8037n^3+9648n^4+4401n^5)\sqrt{3+9n+9n^2}\\
&-(16+60n+57n^2)\sqrt{3n}(-56-555n-1935n^2-1620n^3\\
&+7173n^4+22788n^5+26649n^6+11907n^7)^{1/2}]/\\
&[8(64+672n+2916n^2+6624n^3+8181n^4+4995n^5+999n^6)],\\
x^{01}_2&=[(32+552n+3132n^2+8037n^3+9648n^4+4401n^5)\sqrt{3+9n+9n^2}\\
&+(16+60n+57n^2)\sqrt{3n}(-56-555n-1935n^2-1620n^3\\
&+7173n^4+22788n^5+26649n^6+11907n^7)^{1/2}]/\\
&[8(64+672n+2916n^2+6624n^3+8181n^4+4995n^5+999n^6)].
\end{align*}
The second root lies outside the interval $[\xi_0,\xi_1]$ and,
therefore, could not lead to $\hat{\lambda}.$ More precisely, for
any $n,$ $x^{01}_2(n)>\xi_1(n),$ as one can see by checking the
inequality for $n=1$ and then checking that the derivative of
$x^{01}_2(n)-\xi_1(n)$ (with respect to n) is positive. In order
for $\lambda_0(x_1)$ to be a valid candidate for $\hat{\lambda}$
we shall need to check that
$\lambda_0(x^{01}_1)=\lambda_1(x^{01}_1)<\lambda_2(x^{01}_1)$ or
rule out all other candidates for $\hat{\lambda}.$ We choose the
second route and proceed to analyze the other intersections of
$\lambda_j.$
\par The solutions of the equation $\lambda_1(x)=\lambda_2(x)$ are
\begin{align*}
x^{12}_1&=\frac{3+6n}{8\sqrt{3+9n+9n^2}},\\
x^{12}_2&=\frac{9+18n}{8\sqrt{3+9n+9n^2}}.
\end{align*}
This time both roots lie outside the interval $[\xi_1,\xi_2].$
More precisely, $x^{12}_1<\xi_1$ and $x^{12}_2>\xi_2,$ which is
clear once one recalls the formulae for $\xi_j(n).$ Thus, this
intersection does not produce any candidates for $\hat{\lambda}.$
\par Turning to the last intersection, we find that
$\lambda_2(x)=\lambda_0(x)$ for
\begin{align*}
x^{20}_1&=[(178+1812n+7101n^2+13455n^3+12357n^4+4401n^5)\sqrt{3+9n+9n^2}\\
&-(13+54n+57n^2)\sqrt{3}(689+7847n+39387n^2+113562n^3\\
&+205110n^4+236718n^5+169641n^6+68607n^7+11907n^8)^{1/2}]/\\
&[8(-131-969n-2835n^2-3870n^3-1809n^4+999n^5+999n^6)],\\
x^{20}_2&=[(178+1812n+7101n^2+13455n^3+12357n^4+4401n^5)\sqrt{3+9n+9n^2}\\
&+(13+54n+57n^2)\sqrt{3}(689+7847n+39387n^2+113562n^3\\
&+205110n^4+236718n^5+169641n^6+68607n^7+11907n^8)^{1/2}]/\\
&[8(-131-969n-2835n^2-3870n^3-1809n^4+999n^5+999n^6)].
\end{align*}
A computation shows that the second root lies outside the interval
$[\xi_0,\xi_2]$ and, thus, irrelevant to the computation of
$\hat{\lambda}.$
\par Let us now establish that the first intersection of $\lambda_2$
and $\lambda_0$ does not lead to $\hat{\lambda}$ either.
\begin{claim}
In the notation of this subsection,
$\lambda_2(x^{20}_1)\neq\hat{\lambda}.$
\end{claim}
\begin{proof}
First, we notice that
\begin{align}
x^{12}_1<x^{20}_1<x^{12}_2. \label{l_20_1}
\end{align}
Now, for any $x\in(x^{12}_1,x^{12}_2),$
$\lambda_1(x)<\lambda_2(x)$ since
$\max_x\{\lambda_1(x)\}=\lambda_1(\xi_1)<\lambda_2(\xi_2)=\max_x\{\lambda_2(x)\},$
and $x^{12}_1$ and $x^{12}_2$ (the ordinates of the intersections
of $\lambda_1$ and $\lambda_2$) lie to the left and to the right
of the interval $[\xi_1,\xi_2]$ respectively.  Therefore,
\eqref{l_20_1} implies that
$\lambda_2(x^{20}_1)>\lambda_1(x^{20}_1),$ which means that this
intersection of $\lambda_2$ and $\lambda_0$ lies above $\lambda_1$
and (its ordinate) can not be $\hat{\lambda}.$
\end{proof}
\par Therefore, the only remaining candidate for $\hat{\lambda},$
which is $\lambda_0(x^{01}_1),$ is $\hat{\lambda}:$
\begin{align*}
\hat{\lambda}&=\lambda_0(x^{01}_1)=\frac{17+63n+63n^2}{16+48n+48n^2}\\
&-\frac{1}{16}\left[\frac{(7+33n+33n^2)^2}{(1+3n+3n^2)^2}+4\left(\frac{9(1+2n)}{\sqrt{3+9n+9n^2}}\right.\right.\\
&+\left\{(32+552n+3132n^2+8037n^3+9648n^4+4401n^5)\sqrt{3+9n+9n^2}\right.\\
&-\sqrt{3n}(16+60n+57n^2)(-56-555n-1935n^2-1620n^3+7173n^4\\
&\hphantom{-\sqrt{3n}(16+60n+57n^2)(}\left.\left.+22788n^5+26649n^6+11907n^7)^{1/2}
\vphantom{\sqrt{9n^2}}\right\}\right/\\
&\hphantom{+}\left\{64+672n+2916n^2+6624n^3+8181n^4+4995n^5+999n^6\right\}
\left.\vphantom{\frac{9(1+2n)}{\sqrt{3+9n+9n^2}}}\right)^2
\left.\vphantom{\frac{(7+33n+33n^2)^2}{(1+3n+3n^2)^2}}\right]^{1/2}
\end{align*}
We make two remarks about $\hat{\lambda}.$ First,
$\hat{\lambda}(n)$ is monotonously increasing (as one can see by
checking the positivity of $\hat{\lambda}(n+1)-\hat{\lambda}(n)$),
which implies that $\hat{\lambda}(n)\geq\hat{\lambda}(1)>1/25,$
and, second, $\lim_{n\to\infty}\hat{\lambda}(n)=2/37,$ which is
the minimal curvature of $W(p,q)$ with $p/q=1.$
\par Let us now describe how to compute $\bar{\lambda}$ - the other
quantity needed to estimate $K_{min}(W(p,q)).$ The number
$\bar{\lambda}$ is computed through the following auxiliary
quantities:
$$A=\begin{pmatrix}
a_0&b_2&b_1\\
b_2&a_1&b_0\\
b_1&b_0&a_2
\end{pmatrix}
=\begin{pmatrix} 8 - \frac{9(p + q)^2}{2(p^2 + pq + q^2)}
&-\frac{p^2 + pq + 10q^2}{4(p^2+pq+q^2)}
&-\frac{10p^2 + pq + q^2}{4(p^2+pq+q^2)}\\
-\frac{p^2 + pq + 10q^2}{4(p^2+pq+q^2)} &4-\frac{9p^2}{8(p^2 + pq
+ q^2)}
&-2-\frac{9pq}{8(p^2 + pq + q^2)}\\
-\frac{10p^2+pq+q^2}{4(p^2+pq+q^2)} &-2-\frac{9pq}{8(p^2+pq+q^2)}
&4-\frac{9q^2}{8(p^2+ pq + q^2)}
\end{pmatrix};$$
\begin{align*}
D_1&=a_1a_2-b_0^2-a_1b_1-a_2b_2+b_0b_1+b_2b_0\\
&=\frac{3(113p^2+86pq+113q^2)}{16(p^2+pq+q^2)},\\
D_2&=a_2a_0-b_1^2-a_2b_2-a_0b_0+b_1b_2+b_0b_1\\
&=\frac{9(38p^2+13pq+45q^2)}{16(p^2+pq+q^2)},\\
D_3&=a_0a_1-b_2^2-a_0b_0-a_1b_1+b_1b_0+b_1b_2\\
&=\frac{9(106p^4+141p^3q+192p^2q^2+77pq^3+60q^4)}{32(p^2+pq+q^2)}.\\
\end{align*}
\par It can be shown that
$$
\bar{\lambda}=\min\left\{\mathbf{x}^{t}A\mathbf{x}|x_0,x_1,x_2\geq
0, x_0+x_1+x_2=1 \right\},
$$
where $$A=\mathfrak{R}(\mathbf{0},0)|_{\mathbb{R}^3}=\begin{pmatrix}a_0&b_2&b_1\\
b_2&a_1&b_0\\b_1&b_0&a_2\end{pmatrix}.$$ Using this
characterization of $\bar{\lambda}$ one sees that
$\bar{\lambda}=\frac{\det A}{D_1+D_2+D_3},$ where
$$\begin{pmatrix}D_1\\
D_2\\D_3\end{pmatrix}=\adj(A)\begin{pmatrix}1\\
1\\1\end{pmatrix}$$ under certain conditions on the elements of
$A$ and $D_i.$
\par
\begin{figure}
\begin{center}
\includegraphics{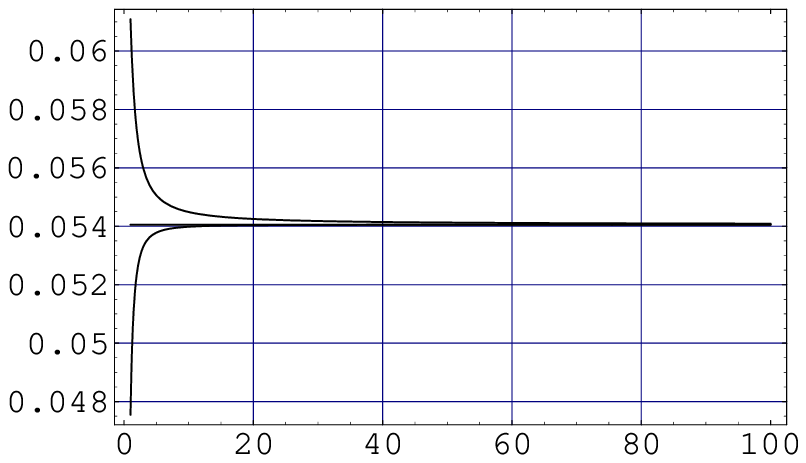}
\end{center}
\caption{$\bar{\lambda}$ approaches $\frac{2}{37}$ from above
while $\hat{\lambda}$ approaches it from below.}
\label{l_hat_v_l_bar}
\end{figure}

 According to proposition 5.11 of \cite{P}, if
$\sum_{j=0}^2(a_j-b_j)\geq 0$ and $D_j>0$ for all $j\in\{0,1,2\},$
$$\bar{\lambda}=\frac{\det A}{D_1+D_2+D_3}.$$
For any $(p,q)\in \mathbb{Z}^{+}\times\mathbb{Z}^{+},$
$\sum_{j=0}^2(a_j-b_j)=\frac{121p^2+85pq+121q^2}{8(p^2+pq+q^2)}>0$
and all $D_j$ are positive (as evident from the above formulae).
Therefore, in this case,
$$\bar{\lambda}=\frac{\det A}{D_1+D_2+D_3}=
\frac{(p^2+pq+q^2)(59p^2-22pq+59q^2)}{772p^4+1127p^3q+1776p^2q^2+977pq^3+676q^4}.$$
In particular, if $(p,q)=(n,n+1)$ (where $n$ is a positive
integer),
$$\bar{\lambda}=\frac{(1+3n+3n^2)(59+96n+96n^2)}{676+3681n+8763n^2+10314n^3+5328n^4}.$$
We remark that the function $\bar{\lambda}(n)$ is decreasing (as
one can see by checking that $\bar{\lambda}'$ is negative), and
$\lim_{n\to\infty}\bar{\lambda}=\frac{3\cdot
96}{5328}=\frac{2}{37},$ which is the minimal curvature of
$W(1,1).$ Figure \ref{l_hat_v_l_bar} compares $\bar{\lambda}$ and
$\hat{\lambda}$ and shows that they both approach the minimal
curvature of $W(1,1),$ but from different sides.
Summarizing this subsection, we can say that the minimal curvature
of the Aloff-Wallach space $W(n,n+1)$ is given by the formula
\begin{align*}
&K_{min}(W(n,n+1))=\frac{17+63n+63n^2}{16+48n+48n^2}\\
&-\frac{1}{16}\left[\frac{(7+33n+33n^2)^2}{(1+3n+3n^2)^2}+4\left(\frac{9(1+2n)}{\sqrt{3+9n+9n^2}}\right.\right.\\
&+\left\{(32+552n+3132n^2+8037n^3+9648n^4+4401n^5)\sqrt{3+9n+9n^2}\right.\\
&-\sqrt{3n}(16+60n+57n^2)(-56-555n-1935n^2-1620n^3+7173n^4\\
&\hphantom{-\sqrt{3n}(16+60n+57n^2)(}\left.\left.+22788n^5+26649n^6+11907n^7)^{1/2}
\vphantom{\sqrt{9n^2}}\right\}\right/\\
&\hphantom{+}\left\{64+672n+2916n^2+6624n^3+8181n^4+4995n^5+999n^6\right\}
\left.\vphantom{\frac{9(1+2n)}{\sqrt{3+9n+9n^2}}}\right)^2
\left.\vphantom{\frac{(7+33n+33n^2)^2}{(1+3n+3n^2)^2}}\right]^{1/2}<\frac{1}{25},
\end{align*}
which proves the left parts of the inequalities from Theorem
\ref{curv_thm}.
\subsection{Maximal Curvature}
\label{ss_max_curv} Let us now describe the computation of the
maximum of the sectional curvature. To do this, we need to compute
numbers $\Lambda_j$ (with $j\in\{0,1,2\}$). In fact, we shall show
that $K_{max}=\Lambda_0.$
\par According to lemma 5.13 of \cite{P}, if $a_1>2d_0-b_0,$
$\Lambda_0=\max\{a_0,a_1,a_2,c_0\}.$ Using the formulae above we
compute (for $W(n,n+1)$):
$$a_1=4-\frac{9n^2}{8(1+3n+3n^2)}>\frac{26+87n+87n^2}{8(1+3n+3n^2)}=2d_0-b_0,$$
which allows us to apply lemma 5.13 (of \cite{P}) and get
\begin{align*}
\Lambda_0&=\max\{a_0,a_1,a_2,c_0\}\\
&=\max\left\{2-\frac{3}{2(1+3n+3n^2)},4-\frac{9n^2}{8(1+3n+3n^2)},\right.\\
&\hphantom{=\max\left\{\right.}\left.4-\frac{9(1+n)^2}{8(1+3n+3n^2)},\frac{3(1+2n)^2}{2(1+3n+3n^2)}\right\}\\
&=4-\frac{9n^2}{8(1+3n+3n^2)}.
\end{align*}
Lemma 5.14 of \cite{P} says that, if $a_1>2d_0-b_0$ (which we
already verified) and
\begin{align*}
\Lambda_0&\geq\max\left\{b_0,b_1,b_2,c_0,\lambda_{max}
\begin{pmatrix}
c_1&\sqrt{2}(\xi_1-\xi_0)\\
\sqrt{2}(\xi_1-\xi_0)&2d_1-b_1&
\end{pmatrix},\right.\\
&\hphantom{\geq\max\left\{b_0,b_1,b_2,c_0,\right.}\left.
\lambda_{max}
\begin{pmatrix}
c_2&\sqrt{2}(\xi_2-\xi_0)\\
\sqrt{2}(\xi_2-\xi_0)&2d_2-b_2&
\end{pmatrix}
\right\},
\end{align*}
(where $\lambda_{max}$ denotes the largest eigenvalue of the
matrix written next to it), $K_{max}=\Lambda_{0}.$ Since we
already know that $\Lambda_0=a_1\geq c_0,$ and all $b_j,$ namely
\begin{align*}
b_0&=-2-\frac{9n(n+1)}{8(1+3n+3n^2)},&
b_1&=-\frac{1+3n+12n^2}{4(1+3n+3n^2)},&
b_2&=-\frac{10+21n+12n^2}{4(1+3n+3n^2)},
\end{align*}
are negative, we only need to check the eigenvalues. Introducing
the notation
\begin{align*}
\nu_1&=\lambda_{max}
\begin{pmatrix}
c_1&\sqrt{2}(\xi_1-\xi_0)\\
\sqrt{2}(\xi_1-\xi_0)&2d_1-b_1
\end{pmatrix},\\
\nu_2&=\lambda_{max}
\begin{pmatrix}
c_2&\sqrt{2}(\xi_2-\xi_0)\\
\sqrt{2}(\xi_2-\xi_0)&2d_2-b_2
\end{pmatrix},
\end{align*}
we compute
\begin{align*}
\nu_1&=\frac{4+12n+33n^2+\sqrt{400+2976n+8640n^2+11664n^3+6561n^4}}{16(1+3n+3n^2)},\\
\nu_2&=\frac{25+54n+33n^2+\sqrt{961+5556n+13014n^2+14580n^3+6561n^4}}{16(1+3n+3n^2).}
\end{align*}
Evidently $\nu_2>\nu_1$; therefore, once we show that
$\Lambda_0=a_1>\nu_2,$ we shall conclude that
$K_{max}=\Lambda_{0}.$ Subtracting $\nu_2$ from $a_1,$ we get
$$
a_1-\nu_2= \frac{39+138n+141n^2-\sqrt{961+5556n+13014n^2+14580n^3+
6561n^4}}{16(1+3n+3n^2)}.
$$
A computation shows that this fraction is always positive.
\par
\begin{figure}
\begin{center}
\includegraphics{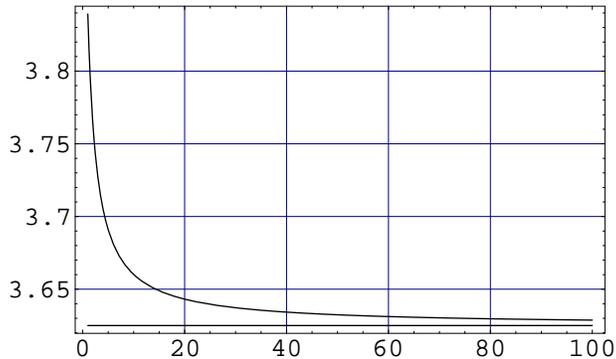}
\end{center}
\caption{$K_{max}(W(n,n+1))$ approaches $K_{max}(W(1,1))$ as
$\frac{n}{n+1}$ tends to one.} \label{L_0}
\end{figure}

 Thus, the maximum curvature of
the Aloff-Wallach space $W(n,n+1)$ is given by the formula
$$K_{max}(W(n,n+1))=4-\frac{9n^2}{8(1+3n+3n^2)},$$ which proves the right
part of the first inequality from Theorem \ref{curv_thm}. We note
that
$$\lim_{n\to\infty}K_{max}(W(n,n+1))=4-\frac{9}{24}=\frac{29}{8}=K_{max}(W(1,1)).$$
Figure \ref{L_0} shows $K_{max}(W(n,n+1))$ together with the
asymptote.

\section{Application: Injectivity Radius Estimates}
\label{ir_sec}
\subsection{Estimating Injectivity Radius from Below}
\label{ir_bel_ss} Cheeger's injectivity radius estimate, [first
obtained in \cite{C}] gives a lower bound on the injectivity
radius in terms of dimension, a lower bound on the volume, a
two-sided bound on the curvature, and an upper bound on the
diameter. We shall use an improved version of this estimate, which
is given in \cite{S} as Theorem IV.3.9(2). For any compact
Riemannian manifold $M,$ let $i(M)$ denote the injectivity radius
of $M,$ $K(M)$ - the set of the sectional curvatures of $M,$ and
$d(M)$ - the diameter of $M.$ Also let $s_{\delta}$ be the
solution of $f''+\delta f=0$ with the initial conditions $f(0)=0$
and $f'(0)=1.$ The above mentioned theorem from \cite{S} says
that, if $K(M)\subset[\delta,\Delta],$
$$i(M)\geq\min\left\{\frac{\pi}{\sqrt{\Delta}},
\pi\frac{Vol(M)}{Vol(S^n)}\cdot
\left[s_{\delta}\left(\min\left\{d(M),
\frac{\pi}{2\sqrt{\delta}}\right\}\right)\right]^{1-n}\right\},$$
where $n$ is the dimension of $M.$ For $\delta>0,$ which the case
for $M=W(p,q),$
$s_{\delta}(t)=\frac{\sin(\sqrt{\delta}t)}{\sqrt{\delta}}.$ Since
$\sin(t)\leq 1,$ $s_{\delta}(t)\leq\frac{1}{\sqrt{\delta}},$ which
implies that $[s_{\delta}(t)]^{1-n}\geq[\sqrt{\delta}]^{n-1}$ for
any $n>1.$ Putting this into the injectivity radius estimate, we
see that, if $K(M)\subset[\delta,\Delta]$ and $\delta>0,$
\begin{equation}
i(M)\geq\min\left\{\frac{\pi}{\sqrt{\Delta}},
\pi\frac{Vol(M)}{Vol(S^n)}\delta^{(n-1)/2}\right\} \label{ch_inj
rad}
\end{equation}
\par Since $$\dim(W(p,q))=\dim(SU(3))-\dim(T(p,q))=8-1=7,$$ when
we apply \eqref{ch_inj rad} to $W(p,q),$ we shall need to put
$$Vol(S^7)=\frac{2\pi^{(7+1)/2}}{\Gamma((7+1)/2)}=\frac{\pi^4}{3}$$
in place of $Vol(S^n)$ and $\delta^3$ in place of
$\delta^{(n-1)/2}.$ Applying \eqref{ch_inj rad} to $M=W(p,q)$ and
using the left part of \eqref{main_est}, we get
\begin{align*}
i(W(p,q))
&\geq\min\left\{\frac{\pi}{\sqrt{\Delta}},\pi\frac{Vol(W(p,q))}{Vol(S^7)}\delta^3\right\}\\
&\geq\min\left\{\frac{\pi}{\sqrt{\Delta}},
\frac{3\sqrt{3}\pi\gcd(p,q)\delta^3}{32\sqrt{p^2+q^2+pq}}\right\}.
\end{align*}
Applying this inequality to $W(n,n+1)$ (where $n$ is a positive
integer) and using the curvature estimates derived in section
\ref{curv_sec}, we get part \ref{ir_n} of Corollary \ref{ir}.
Applying the inequality to $W(1,1),$ where $\delta=2/37$ and
$\Delta=29/8,$ we see that
$$
i(W(1,1))\geq\min\left\{\frac{\pi}{\sqrt{29/8}},
\frac{3\sqrt{3}\pi(2/37)^3}{32\sqrt{1^2+1^2+1\cdot1}}\right\}=\frac{3\pi}{4\cdot37^3}\geq4.65\cdot10^{-5},
$$
which proves item \ref{ir_1} of Corollary \ref{ir}.
\subsection{Estimating Injectivity Radius from Above}
\label{ir_ab_ss} In \cite{B}, Berger proved that
\begin{equation}
Vol(M^m,g)\geq \left[\frac{i(M^m,g)}{\pi}\right]^m\cdot Vol(S^m).
\label{berger}
\end{equation}
This result is known as the Berger isoembolic inequality and its
expository account can be found, for example, in \cite{S} as
Theorem VI.2.1. Since we have an upper bound on the volumes of the
Aloff-Wallach spaces \eqref{main_est}, we can apply \eqref{berger}
to $M^m=W(p,q)$ and get upper bounds on their injectivity radii.
Since (as we saw in the previous subsection) $\dim(W(p,q))=7$ and
$Vol(S^7)=\pi^4/3,$ combining the right part of \eqref{main_est}
with \eqref{berger}, we get
$$i(W(p,q))\leq\left[\frac{3\sqrt{3}\gcd(p,q)}{2\sqrt{p^2+q^2+pq}}\right]^{1/7}\cdot\pi,$$
and Corollary \ref{ir}(\ref{ir_ab}) is proven. In particular, we
see that, if at least one of the indices $p$ and $q$ goes to
infinity, the injectivity radius of $W(p,q)$ tends to zero.
\par Applying this estimate to the family of spaces
considered by Huang in \cite{H}, which are $W(i,i+1)$ with $i$
sufficiently large, we see that the family contains spaces with
arbitrary small injectivity radii, and, thus, we conclude that
there can not be a common lower injectivity radius for this
family, getting an alternative prove of Huang's result mentioned
in the introduction.

\par \textsc{Detartment of Mathematics, Columbia University, New York, NY 10027}
\par \emph{E-mail address}: \texttt{alyurov@math.columbia.edu}
\end{document}